\documentclass[a4paper,12pt]{amsart}
\usepackage[utf8]{inputenc}
\usepackage[english]{babel} 

\usepackage[T1]{fontenc}

\usepackage{amsthm}
\usepackage{amsmath}
\usepackage{amsfonts}
\usepackage{mathtools}
\usepackage[colorlinks]{hyperref}
\usepackage{enumerate}

\def\bdi{\begin{diagram}}
\def\edi{\end{diagram}}

\newtheorem{thm}{Theorem}[section]
\newtheorem{cor}[thm]{Corollary}
\newtheorem{lem}[thm]{Lemma}
\newtheorem{prop}[thm]{Proposition}
\theoremstyle{definition}
\newtheorem{defi}[thm]{Definition}
\newtheorem{defis}[thm]{Definitions}
\newtheorem{conj}[thm]{Conjecture}
\newtheorem{conv}[thm]{Convention}
\newtheorem{nota}[thm]{Notation}
\newtheorem{rem}[thm]{Remark}
\newtheorem{rems}[thm]{Remarks}
\newtheorem{exa}[thm]{Example}
\newtheorem{exas}[thm]{Examples}
\newtheorem{prob}[thm]{Problem}
\newtheorem{probs}[thm]{Problems}
\newtheorem{ques}[thm]{Question}
\newtheorem{sett}[thm]{Setting}
\newtheorem{sit}[thm]{}

\newcommand{\sing}{\operatorname{{\rm sing}}}

\newcommand{\Ker}{ \operatorname{{\rm Ker}}}

\newcommand{\ED}{ \operatorname{{\rm ED}}}

\newcommand{\Aut}{ \operatorname{{\rm Aut}}}

\newcommand{\LND}{\operatorname{{\rm LND}}}

\newcommand{\SL}{ \operatorname{{\rm SL}}}

\def\pr{\mathop{\rm pr}}

\def\reg{{\mathop{\rm reg}}}
\def\sing{{\mathop{\rm sing}}}

\def\codim{\mathop{\rm codim}}

\def\lto{\longrightarrow}

\renewcommand{\epsilon}{\varepsilon}

\def\and{\quad\mbox{and}\quad}

\newcommand{\G}{\ensuremath{\mathbb{G}}}

\newcommand{\bk}{{\ensuremath{\rm \bf k}}}

\newcommand{\hS}{{\hat S}}
\newcommand{\bX}{{\bar X}}

\newcommand{\hZ}{{\hat Z}}

\newcommand{\tQ}{{\tilde Q}}

\newcommand{\bH}{{\bar H}}

\newcommand{\bZ}{{\bar Z}}

\newcommand{\cB}{{\ensuremath{\mathcal{B}}}}

\newcommand{\cG}{{\ensuremath{\mathcal{G}}}}

\newcommand{\cA}{{\ensuremath{\mathcal{A}}}}

\newcommand{\cH}{{\ensuremath{\mathcal{H}}}}

\newcommand{\cN}{{\ensuremath{\mathcal{N}}}}

\newcommand{\id}{{\rm id}}

\renewcommand{\rho}{\varrho}

\def\bals#1\eals{\begin{align*}#1\end{align*}}
\def\bal#1\eal{\begin{align}#1\end{align}}

\def\SAut{\mathop{\rm SAut}}
\def\SL{\mathop{\rm SL}}
\def\Sp{\mathop{\rm Sp}}

\def\A{{\mathbb A}}

\def\PP{{\mathbb P}}

\renewcommand{\phi}{\varphi}

\newcommand{\bnum}{\begin{enumerate}}
\newcommand{\enum}{\end{enumerate}}

\newcommand{\brem}{\begin{rem}}
\newcommand{\brems}{\begin{rems}}
\newcommand{\erem}{\end{rem}}
\newcommand{\erems}{\end{rems}}
\newcommand{\bprob}{\begin{prob}}
\newcommand{\eprob}{\end{prob}}
\newcommand{\bprobs}{\begin{probs}}
\newcommand{\eprobs}{\end{probs}}
\newcommand{\bques}{\begin{ques}}
\newcommand{\eques}{\end{ques}}
\newcommand{\bexa}{\begin{exa}}
\newcommand{\bexas}{\begin{exas}}
\newcommand{\eexa}{\end{exa}}
\newcommand{\eexas}{\end{exas}}
\newcommand{\bdefi}{\begin{defi}}
\newcommand{\edefi}{\end{defi}}
\newcommand{\bdefis}{\begin{defis}}
\newcommand{\edefis}{\end{defis}}
\newcommand{\bcor}{\begin{cor}}
\newcommand{\ecor}{\end{cor}}
\newcommand{\blem}{\begin{lem}}
\newcommand{\elem}{\end{lem}}
\newcommand{\bconv}{\begin{conv}}
\newcommand{\econv}{\end{conv}}
\newcommand{\bconj}{\begin{conj}}
\newcommand{\econj}{\end{conj}}
\newcommand{\bprop}{\begin{prop}}
\newcommand{\eprop}{\end{prop}}
\newcommand{\bthm}{\begin{thm}}
\newcommand{\ethm}{\end{thm}}
\newcommand{\bnota}{\begin{nota}}
\newcommand{\enota}{\end{nota}}
\newcommand{\bsit}{\begin{sit}}
\newcommand{\esit}{\end{sit}}
\newcommand{\be}{\begin{equation}}
\newcommand{\ee}{\end{equation}}
\newcommand{\bproof}{\begin{proof}}
\newcommand{\eproof}{\end{proof}}
\newcommand{\bsett}{\begin{sett}}
\newcommand{\esett}{\end{sett}}
\def\ba{\begin{array}}
\def\ea{\end{array}}

\begin{document}

\title{Embedding theorems for flexible varieties}

\author{Shulim Kaliman}
\address{
University of Miami, Department of Mathematics, Coral Gables, FL 33124, USA}
\email{kaliman@math.miami.edu}

\date{\today}
\maketitle

\begin{abstract}  Let $Z$ be an affine algebraic variety and $X$ be a smooth flexible variety. 
We develop some criteria under which $Z$ admits a closed embedding into $X$.
In particular, we show that if
$\dim X \geq  \max(2\dim Z+1, \dim TZ)$ 
and $X$ is isomorphic (as an algebraic variety)  to a special linear group
or  to a symplectic group, then $Z$ admits a closed embedding into $X$.
\end{abstract}

\thanks{
{\renewcommand{\thefootnote}{} \footnotetext{2020
\textit{Mathematics Subject Classification:}
14E25, 14L30, 14R10. \mbox{\hspace{11pt}}\\{\it Key words}:
closed embedding, injective immersion, affine algebraic variety, flexible variety, special linear group,
symplectic group, linear algebraic group.}}}

\section{Introduction} 

All algebraic varieties which appear in this paper are considered over an algebraically closed field $\bk$ of characteristic zero.
If $Z$ is an affine algebraic variety and $TZ$ is its Zariski tangent bundle then we call $\ED(Z)=\max(2\dim Z+1, \dim TZ)$ the embedding dimension of $Z$.
Holme's theorem \cite[Theorem 7.4]{Hol} (later rediscovered in \cite{Ka91} and \cite{Sr}) states that $Z$ admits a closed embedding
into any affine space $\A^n$ with $n \geq \ED(X)$. In the smooth case (when $\ED(Z)=2\dim Z +1$) this fact was proven earlier by Swan \cite[Theorem 2.1]{Swan}.  The
latter result is sharp -
examples of smooth irreducible $d$-dimensional affine algebraic varieties  with $d \geq \frac{n}{2}$ such that  they do no admit closed embeddings in $\A^{n}$ were constructed
in  \cite{BMS}. Recently Feller and van Santen \cite{FvS21}
proved that if $X$ is an affine variety isomorphic to a simple linear algebraic group and $Z$ is smooth, then $Z$ admits a closed
embedding into $X$, provided that $\dim X>\ED (Z)$.  They also proved  that for every $n$-dimensional algebraic group $G$
(with $n>0$)     
there exist    smooth irreducible $d$-dimensional affine algebraic varieties  with $d \geq \frac{n}{2}$ such that  they do not admit closed embeddings in $G$     \cite[Corollary 4.4]{FvS21}.
In particular, their embedding result is optimal if the dimension of $X$ is even. 
However, they did not know whether their result is sharp in the case  the dimension of $X$ is odd and 
a  specific question posed in \cite{FvS21} asks
whether a smooth affine algebraic variety of dimension 7 can be embedded properly into $\SL_4(\bk)$. 
We  consider a more general situation.  Namely, starting from dimension 2 affine spaces and linear algebraic groups without nontrivial 
characters are examples
of so-called flexible varieties. Recall  that  a normal quasi-affine variety $X$ of dimension at least 2 is flexible 
if $\SAut (X)$ acts transitively on the smooth part
$X_{\rm reg}$ of $X$ where $\SAut (X)$ is
the subgroup of the group $\Aut (X)$ of algebraic automorphisms of $X$ generated by all one-parameter unipotent subgroups  (in what follows  one-parameter unipotent groups will be called $\G_a$-groups and
$\G_a^m$ will stand for the $m$-th power of a $\G_a$-group).
The main results of this paper are the following.

\bthm\label{int.t1} Let $X$ be a smooth flexible variety equipped with a 
$\G_a^m$-action such that the minimal dimension of its orbits is $n$. 
Suppose that $Z$ is an affine variety such that $\dim Z \leq n$
and $\ED (Z)\leq \dim X$. Then there exists a closed embedding of $Z$ into $X$.
\ethm

\bthm\label{int.t2} Let $X$ be isomorphic (as an algebraic variety)
to a connected  linear algebraic group  $G\ne \G_a$ without nontrivial characters.
Suppose that $G'\simeq \G_a^{m'}$ and $G''\simeq \G_a^{m''}$ are subgroups of $G$
such that $G'\cap G''$ coincides with the identity element of $G$.
Let $Z$ be
an affine  algebraic variety such that $\dim Z \leq m'+m''$
and $\ED (Z)\leq \dim X$. Then there exists a closed embedding of $Z$ into $X$.
\ethm

Theorems  \ref{int.t1} and \ref{int.t2} imply the following.

\bcor\label{int.c1}
  Let $X$ be a smooth flexible variety equipped with a free $\G_a^l$-action.
 Let $Z$ be an affine algebraic variety of dimension at most $n+l$ 
 such $\dim X+n \geq \ED (Z)$.
 Suppose that $\psi : X\times \A^n \to Y$ is a finite morphism onto a normal variety $Y$  and 
$S$ is a closed subvariety of $Y$ such that it contains $Y_{\sing}$ and $\dim Z < \codim_Y S$.
Then $Z$ admits 
a closed embedding into $Y$ with the image contained in $Y\setminus S$.  
\ecor

\bcor\label{int.c2} Let $X$ be isomorphic (as an algebraic variety)
either to a special linear group  $\SL_n(\bk)$  or to 
a symplectic group $\Sp_{2n} (\bk)$ and $Z$ be
an affine algebraic variety such that
$\ED (Z)\leq \dim X$. Then there exists a closed embedding of $Z$ into $X$.

\ecor

In particular, the question of Feller and van Santen has  a   positive answer. 
Corollary \ref{int.c2} can be extended to semi-simple Lie groups whose Lie algebras
are direct sums of simple Lie algebras with Dynkin diagrams $A_n$ or $C_n$. In fact, we have more.

\bcor\label{int.c3} Let $Z$ be an affine algebraic variety, $X$ be an algebraic variety of the form 
$\A^{n_0}\times G_{1} \times  G_2 \times \ldots   \times G_l$
where each $G_i$ is
either  $\SL_{n_i}(\bk)$  or  $\Sp_{2n_i} (\bk)$.
Suppose that $\varphi : X \to Y$ is a finite morphism into a normal variety $Y$, $\ED (Z)\leq \dim Y$ and 
$S$ is a closed subvariety of $Y$ containing $ Y_{\sing}$ such that $\dim Z < \codim_Y S$.
Then $Z$ admits 
a closed embedding into $Y$ with the image contained in $Y\setminus S$. 
\ecor

The proofs of Theorems \ref{int.t1} and  \ref{int.t2} are  heavily based on the theory of flexible varieties and
the technique developed in \cite{AFKKZ}, \cite{Ka20}, \cite{KaUd} and \cite{Ka21} 
whose survey can be found in Section 2.
As a part of this survey we describe injective immersions of affine algebraic varieties into smooth flexible varieties.
In section 3 we consider a surjective morphism $\varphi : \A^t \to X$
(every flexible variety $X$ admits such morphism) and 
for a closed subvariety $Z$ of $\A^t$ we
develop a criterion of properness of the morphism $\varphi|_Z : Z \to X$.
Checking the validity of the criterion for injective immersions under the assumptions of
Theorems \ref{int.t1} and \ref{int.t2} 
we prove these theorems in sections 4 and 5.

{\em Acknowlegement.}  The author is grateful to L. Makar-Limanov, Z. Reichstein
 and A. Dvorsky
for  useful consultations and the referee  who
simplified some proofs and caught
mistakes in the original versions of this  paper.

\section{Flexible varieties}

Let us start with the main definitions for the theory of flexible varieties.

\bdefi\label{pre.d1}  (1) Given an irreducible algebraic variety $\cA$ and
a map $\varphi:\cA\to\Aut(X)$ 
we say that $(\cA,\phi)$
is an {\em algebraic family of automorphisms of $X$} if the induced map
$\cA\times X\to X$, $(\alpha,x)\mapsto \varphi(\alpha).x$ is a morphism (see \cite{Ra}).

(2)  If we want to emphasize additionally that  $\varphi (\cA)$ is contained in a subgroup $G$ of $\Aut (X)$, then we say that
$\cA$ is an {\em algebraic  $G$-family} of automorphisms of $X$.            

(3) In the case when $\cA$ is a connected algebraic group and the induced map 
$\cA\times X\to X$ is not only a morphism but also an action of $\cA$ on $X$ we call this family a {\em connected algebraic subgroup} of $\Aut (X)$.

(4)
Following \cite[Definition 1.1]{AFKKZ} we call a subgroup $G$ of $\Aut (X)$ 
{\em algebraically generated} 
if it is generated as an abstract group by a family 
$\cG$ of connected algebraic subgroups of $\Aut (X)$.
\edefi

\bdefi\label{pre.d2}  (1) A nonzero derivation $\delta$ on the ring $A$ of regular functions on an affine algebraic variety $X$ is called 
{\em  locally nilpotent}
if for every $ a \in A$ there exists a natural $n$ for which $\delta^n (a)=0$. 
This derivation can be viewed as a vector field on $X$ which
we also call {\em locally nilpotent}. The set of all locally nilpotent vector fields on $X$ will be denoted by $\LND (X)$. 
The flow of $\delta \in \LND (X)$ 
 is an algebraic $\G_a$-action on $X$, i.e., the action of the group $(\bk, +)$ 
which can be viewed as a one-parameter unipotent group $U$ in the group $\Aut (X)$ of all algebraic automorphisms of $X$.
In fact, every $\G_a$-action is a flow of a locally nilpotent vector field (e.g, see \cite[Proposition 1.28]{Fre}).

(2) If $X$ is a quasi-affine variety, then an algebraic vector field $\delta$ on $X$ is called 
{\em locally nilpotent}
if $\delta$ extends
to a locally nilpotent vector field  $\tilde \delta$ on some affine algebraic variety $Y$  containing $X$ as an open subset such
that $\tilde \delta$ vanishes on $Y\setminus X$ where ${\rm codim}_C (Y \setminus X) \geq 2$.  Note that under this assumption
$\delta$ generates a $\G_a$-action on $X$ and we use again the notation $\LND (X)$ for the set of all locally nilpotent vector
fields on $X$.\edefi

\bdefi\label{pre.d3}

(1) For every locally nilpotent vector fields $\delta$ and each function $f \in \Ker \delta$ from its kernel the field
$f\delta$ is called a {\em replica }
 of $\delta$. Recall that such  a    replica is automatically locally nilpotent.

(2)  Let $\cN$ be a set of locally nilpotent vector fields on $X$ and $G_\cN \subset  \Aut (X)$ denotes the group generated
by all flows of elements of $\cN$. We say that $G_\cN$ {\em is generated by $\cN$}.

(3) A collection of locally nilpotent vector fields $\cN$ is called  {\em saturated}
 if $\cN$ is closed under conjugation by elements in $G_\cN$
and for every $\delta \in \cN$ each replica of $\delta$ is
also contained in $\cN$.

\edefi

\bdefi\label{pre.d4}  Let $X$ be a normal quasi-affine algebraic variety of dimension at least 2,
$\cN$ be a saturated set of locally nilpotent vector fields on $X$ and  $G=G_\cN$ be the group
generated by  $\cN$.
Then $X$ is called $G$-flexible if for every
 point $x$ in the smooth part $X_{\rm reg}$ of $X$ the vector space $T_xX$ is generated 
by  the values of locally nilpotent vector fields from $\cN$ at $x$
(which is equivalent to the fact that $G$ acts transitively on $X_{\rm reg}$ \cite[Theorem 2.12]{FKZ}). In the case of $G=\SAut (X)$ we call $X$ flexible
without referring to $\SAut (X)$ (recall that $\SAut (X)$ is the subgroup of $\Aut X$ generated by all one-parameter unipotent subgroups).
\edefi
\bnota\label{pre.n1} Further in this paper $X$ is always a smooth quasi-affine variety and $G$ is  a group acting transitively on $X$
 such that $G$ is algebraically generated by a collection  $\cG$ of connected algebraic subgroups of $G$.
Given a sequence $\cH=(H_1, \ldots, H_s)$ of elements of $\cG$
we consider the map
\be\label{pre.eq1}
\Phi_\cH : H\times X \lto X\times X, \, (h_s,\ldots,h_1,x)\mapsto
((h_s\cdot\ldots\cdot h_1).x ,x)
\ee
where $H =H_s\times \ldots\times H_1$. 
By $\varphi_\cH : H \lto X$ we denote the restriction of $\Phi_\cH$ to $H \times x_0$ where
$x_0$ is a fixed point of $X$.
\enota

\bprop\label{pre.p1} Suppose that $\cG$ is  closed under  conjugation by $G$.

Then a sequence $\cH=(H_1, \ldots, H_s)$ can be chosen so that for a dense open subset $U$ of $H$
the morphism $\Phi_\cH$ is smooth on $U\times X$ (in particular, $\varphi_\cH$ is smooth on $U$).

{\rm (2)} Let $\cH=(H_1, \ldots, H_s)$ be as in (1)
and $H$ be any element $\cG$. Then
the sequence  $H_1,\ldots, H_m, H$ (resp. $H, H_1,\ldots, H_m$) satisfies  the conclusions of (1) as well.

{\rm (3)} 
Furthermore, increasing the number of elements in $\cH$
one can suppose that the codimension of $H\setminus U$ in $H$ is arbitrarily large.

\eprop

\bproof The first statement follows from \cite[Proposition 1.16]{AFKKZ}, 
the second statement
follows from \cite[Proposition 1.10]{Ka20}) and the third one 
from \cite[p. 778, footnote]{AFKKZ}.
\eproof

We shall use the notion of a perfect (algebraic) $G$-family of automorphisms of $X$  (see \cite[Definition 2.7]{Ka21}).
Without stating the formal definition of such families we need to emphasize some of their properties.

\bprop\label{pre.p2} {\rm (\cite[Proposition 2.8]{Ka21})} Let $\cA$ be a perfect $G$-family of automorphisms of a smooth $G$-flexible variety $X$
and $H_0\in \cG$. Then $H_0\times \cA$ and $\cA\times H_0$ are also  perfect $G$-families of automorphisms of $X$.
Furthermore,  $\cA$ satisfies  the transversality theorem  (\cite[Theorem 1.15]{AFKKZ},  see also \cite[Theorem 2.2]{Ka21}),
 e.g., if $Z$ and $W$ are subvarieties   of $X$ with $\dim Z+\dim W< \dim X$, then
 one has $\alpha (Z) \cap W =\emptyset$ for a general $\alpha \in \cA$.
\eprop

\bthm\label{pre.t1}  
Let $X$ be a smooth quasi-affine $G$-flexible variety, $\cA$ be a perfect $G$-family of automorphisms of $X$,
$Q$ be a normal algebraic variety and $\rho : X \to Q$ be a dominant morphism.
Suppose that  $Q_0$ is a smooth open dense subset of $Q$,  
$X_0$ is an open subset of $X$ contained in $\rho^{-1} (Q_0)$         and
\be\label{pre.eq2} X_0\times_{Q_0} X_0= 2\dim X -\dim Q.\ee
Let $Y$ be the closure of $\bigcup_{x \in X_0} \Ker \{ \rho_* : T_x X_0 \to T_{\rho (x)} Q_0\}$ in $TX$ and
\be\label{pre.eq3}      \dim Y =2\dim X -\dim Q.            \ee
Let $Z$ be a locally closed reduced subvariety of $X$ with $ \ED (Z) \leq \dim Q$ and $\dim Z < \codim_{\rho^{-1}(Q_0)}( \rho^{-1}(Q_0)\setminus X_0)$.
Then for a general element $\alpha \in \cA$ 
 the morphism $\rho|_{\alpha (Z)\cap X_0}:  \alpha (Z) \cap X_0 \to Q_0$ is an injective immersion.
\ethm

\bproof  In the case of $X_0= \rho^{-1} (Q_0)$ the statement  is 
the combination of \cite[Theorem 2.6]{Ka21} and \cite[Proposition 2.8(5)]{Ka21}. 
In the general case the proof goes without change if one observes that $\alpha (Z)$ does not meet    $\rho^{-1}(Q_0)\setminus X_0$ for a general $\alpha \in \cA$
by the transversality theorem.     \eproof

\bprop\label{pre.p3} Let the assumptions and conclusions of Proposition \ref{pre.p1} hold. Suppose that $H$ itself is
an $F$-flexible variety.
Let $Z$ be a locally closed reduced subvariety of $H$ with $ \ED (Z) \leq \dim X$ (and by the conclusions of Proposition \ref{pre.p1}  with $\dim Z < \codim_H (H \setminus U)$).
Then 
for a general element $\beta \in \cB$ in any perfect $F$-family $\cB$ of automorphisms of $H$
 the morphism $\varphi_\cH|_{\beta (Z)}:  \beta (Z)  \to X$ is an injective immersion.
\eprop

\bproof  
Since $\varphi_\cH|_U: U \to X$ is a smooth morphism Formulas \eqref{pre.eq2} and \eqref{pre.eq3} hold with $\rho : X\to Q, Q_0$ and $X_0$ 
replaced by $\varphi_\cH : H \to X, X$ and $U$,
respectively.
Hence, the desired conclusion follows from Theorem \ref{pre.t1}.
\eproof

\bcor\label{pre.c1} Let the assumptions and conclusions of Proposition \ref{pre.p1} hold
and $Z$ be an affine algebraic variety with $ \ED (Z) \leq \dim X$ (and by the conclusions of Proposition \ref{pre.p1}  with $\dim Z < \codim_H (H \setminus U)$). 
Suppose that each element of $\cG$
is a unipotent group, i.e. $H \simeq \A^t$ where $t \geq \dim X$. Then $Z$ can be treated as a closed subvariety of $H$
and for a general element $\beta \in \cB$ in any perfect $F$-family $\cB$ of automorphisms of $H$
 the morphism $\varphi_\cH|_{\beta (Z)}:  \beta (Z)  \to X$ is an injective immersion.
\ecor

\bproof The first statement follows from Holme's theorem
and the second from Proposition \ref{pre.p3}.
\eproof

Since every smooth flexible variety $X$ admits a morphism  $\varphi_\cH: H   \to X$ as  in Corollary \ref{pre.c1} we have the following. 

\bthm\label{pre.t2} {\rm (\cite[Theorem 3.7]{Ka21})}  Let $Z$ be an affine algebraic variety and $X$ be a smooth quasi-affine flexible variety of dimension at least $\ED (Z)$.
Then $Z$ admits an injective immersion into $X$.
\ethm

\brem\label{pre.r1} It is worth mentioning that if $\varphi : Z \to X$ is an  injective immersion, then it may happen that $Z$ is not isomorphic to $\varphi (Z)$.
As an example one can consider the morphism $\A^1\setminus \{ 1\} \to \A^2,\, t\mapsto (t^2-1, t(t^2-1)) $.
It maps
$\A^1\setminus \{ 1\}$ onto  the polynomial curve given  in $\A^2$ by the equation $y^2=x^2(x+1)$.                  

\erem

We have also in our disposal the following slightly improved version of \rm (\cite[Theorem 3.2]{Ka21}.

\bthm\label{pre.t3} 
Let $\psi : X \to Y$ be a finite morphism where $X$ is a smooth  flexible variety and $Y$ is normal. Let $Z$ be a quasi-affine algebraic variety
which admits a closed embedding in $X$  and has $\ED (Z)\leq \dim X$. Suppose also that $S$ is a closed subvariety of $Y$ such that it
contains $Y_{\sing}$ and 
 $\dim Z < \codim_Y S$. Then $Z$ admits 
a closed embedding in $Y$ with the image contained in $Y\setminus S$. 
\ethm

\bproof 
One can treat $Z$ as a closed subvariety of $X$. By  \cite[Theorem 1.15]{AFKKZ}   there exists an algebraic family $\cA$
of automorphisms of $X$ such that for a general $\alpha\in \cA$ the variety $\alpha (Z)$ does not meet $\psi^{-1} (S)$.
By Proposition \ref{pre.p2} enlarging $\cA$ we can suppose that it is a perfect family.
Theorem \ref{pre.t1} and \cite[Proposition 2.9]{Ka21}  imply now
that $\psi|_{\alpha (Z)} : \alpha (Z) \to Y_\reg \subset Y$ is an injective immersion. Since $\psi$ is finite $\psi|_{\alpha (Z)}$ is also proper.  Hence, we are done.
\eproof

\section{Criterion of Properness}

\bnota\label{cri.n1}
In this section an affine space $H=\A^t$ is equipped with a fixed coordinate system.
This coordinate system defines an embedding $H\hookrightarrow \PP^t=\bH$
and we let $D=\bH\setminus H$. By $\varphi : H \to X$ 
we denote a surjective morphism onto a smooth quasi-affine algebraic variety $X$
(of positive dimension)
with irreducible fibers and by  $\psi : \bH \dashrightarrow \bX$ we denote the rational
map into a completion $\bX$ of $X$ extending $\varphi$.
\enota

\bprop\label{cri.p1} 
Let $\pi : Y\to \bH$ be a resolution of the indeterminacy set of $\psi$,
 (i.e.,  $H$ is naturally contained as an open dense subset in $Y$ and $\chi:=\psi \circ \pi : Y\to \bX$ is a proper morphism). 
 Let $V=\chi^{-1} (X) \setminus H$
 and  $W= \pi (V )$.
 Suppose that $Z$ is a closed subvariety of $H$ and $\bZ$ is its closure
in $\bH$. Then  $\varphi|_Z : Z \to X$ is a proper morphism if and only if $\bZ\cap W=\emptyset$.
\eprop

\bproof 
Let  $\hZ=\pi^{-1} (\bZ) \cap V$. 
Note that $\varphi|_Z=\chi|_Z$ is proper  if and only if $\hZ=\emptyset$.
Note also that $\pi (\hZ) =\bZ\cap W$.
In particular,
$\hZ=\emptyset$ if and only if  $\bZ\cap W=\emptyset$.
This yields  
the desired conclusion.
\eproof

\bdefi\label{cri.d1} We call the set $W$ as in Proposition \ref{cri.p1}
{\em the improperness set }of $\varphi$.
\edefi

It is easy to see that if $\dim Z > \codim_DW$, then  $\bZ\cap W\ne \emptyset$.
Hence, in the rest of this section we describe some conditions which guarantee
that $ \codim_DW$ is sufficiently large.

\bprop\label{cri.p2} Let Notation \ref{cri.n1} hold and $G$
be a subgroup of the group  of affine transformations of $H$
(in particular, the natural action of $G$ extends to $\bH$).
Suppose that $G$ acts on $X$ so that the morphism
$\varphi : H \to X$ is equivariant. Then $\bX$ and a resolution 
 $\pi : Y\to \bH$ of the indeterminacy points of $\psi$
 can be chosen such that $G$ acts on $Y$ and
 $\pi$ is equivariant.
\eprop

\bproof By Sumihiro's theorem \cite{Su} we can suppose that
the $G$-action on $X$ extends to a $G$-action on $\bX$.
Then $\psi$ is an equivariant rational map into a complete variety
and the desired conclusion follows from the   Reichstein-Youssin theorem \cite{ReYo}.
\eproof

\bprop\label{cri.p3}  Under the assumptions of Proposition \ref{cri.p2} suppose that 
$G$ acts on $H$ by translations (in particular, the $G$-action on $D$ is trivial) and
the minimal dimension of orbits of $G$ in $X$ is $m$.
Then the codimension of the  improperness set  $W$ of $\varphi$ in $D$ is at least $m$.
\eprop

\bproof Let $U$ be an irreducible component of $V$ where $V$ is as in Proposition \ref{cri.p1}.
Since $\chi|_U : U \to X$ is equivariant
the dimension of a general $G$-orbit in $U$ is at least  $m$.
Since the $G$-action on $D$ is trivial a general fiber of $\pi|_U : U \to \pi (U) \subset D$
contains a $G$-orbit. Hence $\dim \pi (U) \leq \dim U-m$.
Since  $\dim U \leq \dim D$ we have 
the desired conclusion.
\eproof

\bprop\label{cri.p4} 
Suppose that the assumptions of Proposition \ref{cri.p2} hold,
$G$ acts on $H$ by translations and
the dimension of  general orbits of $G$ in $X$ is $n$.  Let $R\subset X$ be the union
of non-general orbits of $G$. Suppose that  $\chi (U)$ is not contained in $R$
for every irreducible component $U$ of $V$ where $V$ is as in Proposition \ref{cri.p1}.
Then the codimension of the  improperness set  $W$ of $\varphi$ in $D$ is at least $n$.
\eprop

\bproof Since $\chi|_U : U \to \chi (U) \subset X$ is equivariant
the dimension of a general $G$-orbit in $U$ is at least the same as the dimension
of general $G$-orbits in $\chi (U)$. By the assumption, the latter dimension is $n$.
Since a general fiber of $\pi|_U : U \to \pi (U) \subset D$
contains a general $G$-orbit one has $\dim \pi (U) \leq \dim U-n\leq \dim D-n$ which concludes the proof.
\eproof

\section{Main Theorem I}

The aim of this section is the following.

\bthm\label{main.t1} Let $X$ be a smooth flexible variety equipped with a 
$\G_a^m$-action such that the minimal dimension of its orbits is $n$. 
Suppose that $Z$ is an affine variety such that $\dim Z \leq n$
and $\ED (Z)\leq \dim X$. Then there exists a closed embedding of $Z$ into $X$.
\ethm

Let us start with the following.

\blem\label{main.l1} Let $G'$ be a $\G_a^{m}$-subgroup of $\SAut (X)$ acting
on $X$. Consider the natural $G'$-action on $X\times X$ given by $(g,x_1,x_2)\mapsto (g.x_1,x_2)$.
Let $\Phi_\cH : H\times X \to H\times X, \, (h,x) \mapsto (h.x,x)$ be as in Proposition \ref{pre.p1}.
Then $\cH$ can be chosen such that $H$ is an affine space equipped with 
a free $G'$-action for which $\Phi_\cH$
  is $G'$-equivariant
(where $G'$ acts on $H\times X$ by $(g,h,x)$ $\mapsto (g.h,x)$).
Furthermore, $H$ can be equipped with a coordinate system such that
$G'$ acts  on $H$ by translations.
\elem

\bproof We can suppose that $\cG$ in Notation \ref{pre.n1} is  the collection of all $\G_a$-subgroups
of $\SAut (X)$ which implies that $H$ is an affine space.
By   Proposition \ref{pre.p1}(2) we can also suppose that
$$\cH=(H_1, \ldots, H_s, H_{s+1}, \ldots, H_{s+m} )$$ where  $H_{s+1}, \ldots, H_{s+m}$
 are commuting $\G_a$-groups generating  $G'$. Let   
  $g'=( h_{s+m}^0, \ldots,  h_{s+1}^0 )\in G'=H_{s+m}\times \ldots \times H_{s+1}$   and
  $h=(h_{s+m}, \ldots,  h_{1} )\in H=H_{s+m}\times \ldots \times H_{1}$.
Suppose that the $G'$-action on $H$ is given 
 by  
 \be\label{main.eq1} (g',h)\mapsto (h_{s+m} h_{s+m}^0, \ldots,  h_{s+1}h_{s+1}^0, h_s, \ldots, h_1).\ee
 Commutativity and Formula \eqref{pre.eq1} imply that $\Phi_\cH (g'.h, x)= (g'. (h.x), x)$ which yields the first statement.
One can equip each $H_i \simeq \A^1$ with a coordinate $\zeta_i$
(with the zero element of $H_i$ corresponding to $\zeta_i=0$). 
This yields the coordinate system  $(\zeta_{s+m}, \ldots, \zeta_{1})$ on $H$.
In this coordinate system the action of $g'$ given by Formula \eqref{main.eq1} is a translation and we are done.
\eproof

\bproof[Proof of Theorem \ref{main.t1}] Let the conclusions of Lemma \ref{main.l1} hold,
$\varphi_\cH : H \to X$  be the restriction of $\Phi_\cH$ to $H\times x_0, \, x_0 \in X$ 
and $U$ be as in Proposition \ref{pre.p1}.
By Holme's theorem we can treat $Z$
as a closed subvariety of $H$ and by
 Proposition \ref{pre.p1}(3) we can suppose $\dim Z< \codim_H (H\setminus U)$. 
By Proposition \ref{cri.p3} and Lemma \ref{main.l1}
 the improperness set $W$ of $\varphi_\cH$ is of codimension at least $n$
in $D=\bH \setminus H=\PP^t\setminus \A^t$.  For any perfect family $\cA$ of 
automorphisms on $H$ and a general 
$\alpha \in \cA$  the morphism
$\varphi_\cH|_{ \alpha (Z)} : \alpha (Z) \to X$ is  an injective immersion by Corollary \ref{pre.c1}.
Let $K=\SL_{s+m}(\bk)$ where  $t=s+m$. Then we have the natural $K$-action on $\bH$ such that $D$ is
invariant under it and the restriction of the action to $D$ is transitive. By Proposition \ref{pre.p2} 
$K\times \cA$ is still a perfect $\SAut(H)$-family of automorphisms of $H$. 
That is, for  a general $(\beta, \alpha) \in K\times \cA$  the morphism
$\varphi_\cH|_{\beta \circ \alpha (Z)} : \beta \circ \alpha (Z) \to X$ is still an injective immersion.
 Let $P$ be the intersection of $D$ with the closure
of $\beta\circ \alpha (Z)$ in $\bH$, i.e., $\dim P \leq n-1$. 
Since the  restriction of the $K$-action to $D$ is transitive, $P$ does not meet $W$  for  a
general $(\beta, \alpha) \in K\times \cA$ by \cite[Theorem 1.15]{AFKKZ}. 
Hence, $\varphi_\cH|_{\beta \circ \alpha (Z)} : \beta \circ \alpha (Z) \to X$ is proper 
by Proposition \ref{cri.p1} and we are done.
\eproof

\bcor\label{main.c1}  Let $X$ be a smooth flexible variety equipped with a free $\G_a^l$-action.
 Let $Z$ be an affine algebraic variety of dimension at most $n+l$ 
 such $\dim X+n \geq \ED (Z)$.
Suppose that $\psi : X\times \A^n \to Y$ is a finite morphism onto a normal variety $Y$  and 
$S$ is a closed subvariety of $Y$ such that it contains $Y_{\sing}$ and $\dim Z < \codim_Y S$.
Then $Z$ admits 
a closed embedding into $Y$ with the image contained in $Y\setminus S$.  
\ecor

\bproof Since $ X\times \A^n$ admits a free $\G_a^{n+l}$-action,  by Theorem \ref{main.t1}
there is a closed embedding of
$Z$ into $ X\times \A^n$.
 Hence, the desired conclusion follows from
Theorem \ref{pre.t3}.
\eproof

\bcor\label{main.c2} Let $X$  be isomorphic (as an algebraic variety)
to a special linear group  $\SL_n(\bk)$ and $Z$ be an affine variety with $\ED (Z) \leq \dim X$.
Suppose also that $\dim Z \leq m=\frac{n^2}{4}$ if $n$ is even and $\dim Z \leq m=\frac{n^2-1}{4}$
if $n$ is odd. Then $Z$ admits a closed embedding into $X$.
\ecor

\bproof  Let $I$ be the identity matrix in $\SL_n (\bk)$. For even $n$ consider the set $G'$ of all matrices of
the form $I +A$ where $A=[a_{ij}]$ is the matrix such that $a_{ij} =0$ as soon as
$i\leq \frac{n}{2}$ or $j> \frac{n}{2}$. 
If $n$ is odd, then we require that $a_{ij} =0$ as soon as
$i\leq \frac{n-1}{2}$ or $j> \frac{n-1}{2}$.  In both cases $G'$ is a $\G_a^m$-group acting freely on $X$
with multiplication given by  $(I+A) \cdot (I+A')=I+ (A+A')$. 
Thus, the desired conclusion follows from Theorem \ref{main.t1}.
\eproof

\section{Main Theorem II}

\bnota\label{main2.n1} In this section $X$ is always isomorphic (as an algebraic variety)
to a connected linear algebraic group  $G\ne \G_a$ without nontrivial characters.
By $\cG$ we denote the collection of all $\G_a$-subgroups of $G$ (the absence
of nontrivial characters implies that such subgroups generate $G$). In particular,
if $\cH=(H_1, \ldots, H_{s})$ is a sequence in $\cG$, then 
the affine space  $H =H_{s}\times \ldots\times  H_{1}$
is equipped with a natural coordinate system as in Lemma \ref{main.l1}.
Recall that we have a morphism  $\Phi_\cH: H\times X \to X\times X$
given by $ \Phi_\cH  (h,x)=
((h_s\cdot\ldots\cdot h_1).x ,x)
$ for  $h =(h_s, \ldots, h_1) \in H_s\times \ldots \times H_1$.
Since we suppose that $G$ acts on $X$  naturally (i.e., $g.x $ coincides with the product $gx$)
$ \Phi_\cH  (h,x)=(hx,x)$ where  $h$ in the right-hand side is treated as the element $h_s \cdot \ldots \cdot h_1$ of $G$.
We also suppose that $G'$ is a $\G_a^{m'}$-subgroup of $G$ which acts on $H$
in the manner described in Lemma \ref{main.l1}. \enota

Our aim is  to strengthen Theorem \ref{main.t1} for such $X$ and, in particular,
to improve Corollary \ref{main.c2}.  
Let us start with some technical facts.

\blem\label{main2.l1}  Let Notation \ref{main2.n1}  hold, $\pr_1 : X\times X\to X$
be the natural projection to the first factor 
and $\Phi_\cH^1=\pr_1\circ \Phi_\cH: H\times X \to X$.
Let $\Lambda : G'\times G\times H \times X\to H\times X, \,(g',g,h,x) \mapsto (g'.h, xg^{-1})$,
$\Delta : G'\times G\times X \times X \to X\times X,  \,
(g',g,x_1,x_2) \mapsto (g'x_1g^{-1}, x_2g^{-1})$ and
$\Delta_1: G'\times G\times X \to X,  \,
(g',g,x) \mapsto g'xg^{-1}$  
be the  $G'\times G$-actions on $H\times X$, $X\times X$ and $X$.
Then $\Phi_\cH$  and $\Phi_\cH^1$ are $G'\times G$-equivariant.
\elem

\bproof Formula \eqref{main.eq1} implies that $ \Phi_\cH  (g'.h,x)=(g'hx,x)$.
Hence, we have 
$$\begin{array}{c} \Phi_\cH (\Lambda (g',g,h,x))=\Phi_\cH (g'.h, xg^{-1}) =(g'hxg^{-1},xg^{-1})\\
=\Delta (g',g, \Phi_\cH (h,x)).\end{array} $$
Thus, $\Phi_\cH$ is equivariant.  Since the morphism $\pr_1$ is also equivariant  we have the desired conclusion.
\eproof

Let $\bX$ be a $\Delta_1$-equivariant completion of $X$
(which implies that $\bX\times \bX$ is a $\Delta$-equivariant completion of $X\times X$).
Then  the proof of Proposition \ref{cri.p2} implies the following.

\blem\label{main2.l2} Let the assumptions of Lemma \ref{main2.l1} hold,
$\bH$ be a   $G'\times G$-equivariant completion of $H\times X$
and  $\Psi: \bH \dashrightarrow \bX\times \bX$ (resp. $\Psi_1: \bH \dashrightarrow \bX$)
be the rational extension of $\Phi_\cH$ (resp. $\Phi_\cH^1$).
Then a resolution $\pi : Y\to \bH$ of the indeterminacy points of $\Psi$ 
can be chosen such that
the $G'\times G$-action on $H\times X$ extends to $Y$ and the morphisms
$\lambda=\Psi \circ \pi:  Y\to \bX\times \bX$ and $\chi=\Psi_1\circ \pi:  Y\to \bX$ are 
$G'\times G$-equivariant.

\elem

\bnota\label{main2.n2}
 From now on we suppose that the conclusions of Lemma \ref{main2.l2} hold
 and we denote the extension of
the $\Lambda$-action on $H\times X$ to $Y$  by the same letter  $\Lambda$ and
the extension of the $\Delta_1$-action to $\bX$ by
the same letter  $\Delta_1$.
For a $\G_a^{m''}$-subgroup $G''$ of $G$ we consider the quotient morphism 
 $\gamma: G\to Q={G''}\backslash G$.
The fiber of this morphism over a  point $q \in Q$ is a right coset of $G''$ denoted by $C_q$.
Fixing an isomorphism $G\simeq X$ we treat $C_q$ as a subset of $X$ and let $H_q=H\times C_q$.
Finally, by $Y_q$ 
we denote the closure of $H_q$ in $Y$.
\enota

\blem\label{main2.l4}  Let Notation \ref{main2.n2} hold 
and $\chi_q : Y_q \to \bX$ be the restriction of $\chi$.
Suppose that $V_q=\chi_q^{-1} (X) \setminus H_q$
and $R$ is a proper closed subvariety of $X$. Then for a general $q \in Q$
there is no irreducible component $U_q$ of $V_q$ with $\chi_q (U_q)$ contained in $R$.
\elem

\bproof 
Note that $V_q=(\chi^{-1} (X) \cap Y_q) \setminus H_q= (\chi^{-1} (X) \setminus (H\times X))  \cap Y_q
=Y_q\cap V $ where $V=\chi^{-1} (X) \setminus (H\times X)$. 
Since  $\bX \setminus X$ is $\Delta_1$-invariant
$\chi^{-1} (\bX \setminus X)$ is 
$\Lambda$-invariant. Since $H\times X$ is also $\Lambda$-invariant,
so is $V=Y\setminus (\chi^{-1} (\bX \setminus X) \cup (H\times X))$.
Note that the $\Lambda$-action yields a transitive action on the collection $\{ H\times C_q\}_{q\in Q}$ and, therefore,
on  $\{ Y_q\}_{q\in Q}$ and, consequently,  on $\{ V_q\}_{q\in Q}$. 
 Thus, $V=\bigcup_{q \in Q} V_q$ is a $\Lambda$-orbit 
of $V_{q_0}$ where $q_0$ is any point in $Q$. Let $q_0$ be the coset $G''$.
Note that
the action of any element of the subgroup $G'\times G''\subset G'\times G$ preserves
$H\times C_{q_0}$ and, therefore, $V_{q_0}$.
Hence, the image of $V_{q_0}$ under the action of $(g',g) \in G'\times G$ is completely determined
by $\tilde \gamma (g)$ where $\tilde  \gamma : G\to G/G''=:\tQ$ is the quotient morphism.
Let $q$ be the
image of $\tilde \gamma (g)$ under the map $\tQ\to Q$ induced $G\to G, \, g \mapsto g^{-1}$. The description of the $\Lambda$-action in Lemma \ref{main2.l1} 
implies that $(g',g).V_{q_0}=V_q$. Note also that every irreducible component $U_{q_0}$ of $V_{q_0}$
is preserved by the action of $G'\times G''$ since the latter subgroup is connected. Hence,
 $(g',g).U_{q_0}$ is a well-defined irreducible component $U_q$ of $V_q$ depending only on $\tilde \gamma (g)$.
 This implies that $\bigcup_{q \in Q} U_{q}$ is the $\Lambda$-orbit of $U_{q_0}$.
Thus, $\chi(\bigcup_{q \in Q} U_{q}) =X$ because $\chi$ is equivariant
and  the $\Delta_1$-action is transitive on $X$.
In particular, $\chi_q (U_q )$ is not contained in $R$ for a general $q\in Q$.
This yields the desired conclusion. \eproof

\blem\label{main2.l5} Let the assumptions of Lemma \ref{main2.l4} hold,
 $q$ be a general point of $Q$
and $C_q=G''g_0$.
Then
$H_q$ is an affine space equipped with a coordinate system such that
in this system the group $G'\times (g_0^{-1}G''g_0)$ acts
on $H_q$ freely by translations.
\elem

\bproof The space $H_q$ is affine since it is isomorphic to $H\times G''$.
Lemma \ref{main.l1} yields a free action of $G'$ on the first factor,
while $g_0^{-1}G''g_0$ acts on the second by multiplications from the right. 
Note also that if $H$ is equipped with a coordinate system from Lemma \ref{main.l1}
and $G''$ with a coordinate system induced by the structure of a $\G_a^{m''}$-subgroup,
then $G'\times g_0^{-1}G''g_0$ acts
on $H_q$ by translations.
Hence, we are done.
\eproof

\blem\label{main2.l6} A completion $\bH$ of $H\times X$ in Lemma \ref{main2.l2} can be chosen such
that for every $q \in Q$ the closure $\bH_q$ of $H_q$ in $\bH$ is a projective space
 that is the completion of $H_q$
associated with the coordinate system from Lemma \ref{main2.l5}.
\elem

\bproof   By \cite[Theorem 3]{Gro58}
the quotient morphism $\gamma: G\to Q$ is a principal $G''$-bundle 
which is locally trivial in the Zariski topology.
Let $\{ Q_i\}$  be a cover of $Q$ by open subsets over which
 $\gamma$ admits sections $\sigma_i : Q_i \to G$.
 The coordinate system on $H$ (from Lemma \ref{main.l1}) allows us to treat $H$ as $\G_a^s$-group. 
Thus,  $\tau : H\times G \to Q$  is a principal $H\times G''$-bundle 
 whose fiber $\tau^{-1} (q)=H_q$ and we have the  trivialization isomorphisms
$$\eta_{i} : Q_i\times H\times G'' \to \tau^{-1} (Q_i), \, (q,h,g'') \mapsto (h, g''\sigma_i(q)) \in H_q$$ 
with the transition functions     $$\kappa_{ij} : Q_{ij}\times H\times G'' \to 
Q_{ij}\times H\times G'', \, (q,h,g'') \mapsto (q,h, g''\sigma_i(q)\sigma_j(q)^{-1}).$$   
Consider  the $G$-action
 on $Q$ such that $g\in G$ sends  $q=G''g_0$ to $G''g_0g^{-1}$ and the 
 set $$S_{ij} =\{ (g',g,q,h,g'') \in G'\times G \times Q_i\times H\times G''| \,  g.q\in Q_j\}.$$
 Then $ \eta_j^{-1} \circ \Lambda \circ (\id, \eta_i) : S_{ij} \to Q_j\times H\times G'' $ is given by
\be\label{main2.eq1} (g',g, q,h,g'') \mapsto  \eta_{j}^{-1}  ((g',g). \eta_{i}(q,h,g''))= (g.q, g'h, g'' \tilde g_{ij}''),\ee  
 where $G''\ni \tilde g_{ij}'' =\sigma_i (q) g^{-1} (\sigma_j (g.q))^{-1}$.
  Equip $H\times G''\simeq \A^t$ (where $t=s+m''$) with the coordinate system $\bar \zeta =(\zeta_1, \ldots, \zeta_t)$
from Lemma \ref{main2.l5}. 
If $\bar \zeta \in \A^t$ are the coordinates of $(h,g'')$ and $\bar \zeta^0(g,q)$ are the coordinates
of $(\bar 0, \tilde g_{ij}'')\in H\times G''$, then the coordinate form of Formula \eqref{main2.eq1} is
\be\label{main2.eq3} (g',g, q,\bar \zeta) \mapsto \eta_{j}^{-1}  ((g',g). \eta_{i}(q,\bar \zeta))= (g.q, \bar \zeta+ \bar \zeta^0(g,q)).\ee 
There is the natural embedding $\A^t \hookrightarrow \PP^t$ where $\PP^t$ is equipped
with the coordinate system $\bar \xi=(\xi_0:\xi_1: \ldots : \xi_t)$ such that $\xi_i=\zeta_i\xi_0$ for $i \geq 1$ and 
$\xi_0\ne 0$.
Since $\kappa_{ij}$ are translations over $Q_{ij}$   the isomorphisms $\eta_{ij}$ extend 
to the trivialization
isomorphisms $\hat \eta_{i} : Q_i\times \PP^t \to\hat \tau^{-1} (Q_i)$ where
$\hat\tau : \widehat{H\times G} \to Q$ is the proectivization of the bundle $\tau : H\times G \to Q$.
For $\hS_{ij}=\{(g',g, q, \bar \xi)\in G'\times G\times Q_i\times \PP^t|\,  g.q\in Q_j\}$ formula \eqref{main2.eq3} admits the extension to the morphism $\hS_{ij}\to  Q_j \times \PP^t$ 
sending $ ((g',g, q,\bar \xi) $ to $(g.q, \bar \xi+ \bar \xi^0(g,q))$ 
where $\bar \xi^0(g,q)=(\xi_0: \xi_1(g,q):\ldots \xi_t(g,q))$ with $\xi_i(g,q)=\zeta_i(g,q)\xi_0$ 
for $i \geq 1$.
Such morphisms yield the morphisms $(\id, \hat \eta_i) (\hS_{ij}) \to \hat \tau_j^{-1} (Q_j)$ which are in turn the
extensions of $\Lambda$ restricted to $(\id, \eta_i)(S_{ij})$. Hence, we have
a $(G'\times G)$-action on $ \widehat{H\times G}$ extending $\Lambda$.  Thus,
 a $(G'\times G)$-equivariant  completion of $ \widehat{H\times G}$  yields $\bH$ which concludes the proof.
\eproof

\bthm\label{main2.t1} Let $X$ be isomorphic (as an algebraic variety)
to a connected linear algebraic group  $G\ne \G_a$ without nontrivial characters.
Suppose that $G'\simeq \G_a^{m'}$ and $G''\simeq \G_a^{m''}$ are subgroups of $G$
such that $G'\cap G''$ coincides with the identity element of $G$.
Let $Z$ be
an affine variety such that $\dim Z \leq m'+m''$
and $\ED (Z)\leq \dim X$. Then there exists a closed embedding of $Z$ into $X$.
\ethm

\bproof 
 Let $q \in Q$, $C_q=G''g_0$, $H_q$ and $Y_q$  be as in Notation  \ref{main2.n2}
 and Lemma \ref{main2.l5} (i.e., $H_q\simeq \A^t$ is an  affine space).
 Consider the group $F = G'\times (g_0^{-1}G''g_0)$  and
 the $F$-actions on $H_q$ and $X$ that are the restrictions
 of $\Lambda$ and $\Delta_1$ from Lemma \ref{main2.l1}, respectively. By Lemma  \ref{main2.l1}
 the morphism $\varphi_q=\Phi_\cH^1|_{H_q}: H_q\to X$ is $F$-equivariant.
 By Lemma \ref{main2.l5} 
 $H_q$ is equipped with a coordinate system
 such that  $F$ acts on $H_q$ by translations.
 Let $\psi_q : \bH_q \dashrightarrow \bX$ be the rational extension of $\varphi_q$
 to the projective space $\bH_q\simeq \PP^t$ which is the completion of $H_q$ associated with the coordinate system.
 By Lemmas   \ref{main2.l2} and \ref{main2.l6} we can suppose that $\pi_q=\pi|_{Y_q} : Y_q \to \bH_q$
 is a $F$-equivariant resolution of the indeterminacy points of $\psi_q$.
 Hence, by Proposition \ref{cri.p4} and Lemma   \ref{main2.l4} 
 we can suppose that the codimension of the improperness set  $W_q$ of $\varphi_q$ in $D_q=\bH_q\setminus H_q$
 is at least the dimension of general orbits of $F$ in $X$.
 Treating $g_0$ as a point in $X\simeq G$ we see that the $F$-orbit of $g_0$
has dimension $m'+m''$.  Thus, the dimension of general $F$-orbits is at least $m'+m''$
and $\codim_{D_q} W_q\geq m'+m''$. 

Let $K=\SL_{t}(\bk)$ and $\cA$ be a perfect family $\cA$ of 
automorphisms on $H_q$. By Holme's theorem we can treat $Z$
as a closed subvariety of $H_q$.
Arguing as in the proof of Theorem \ref{main.t1} we see that for
a general  $(\beta, \alpha) \in K\times \cA$  the morphism
$\varphi_q|_{\beta \circ \alpha (Z)} : \beta \circ \alpha (Z) \to X$ is an injective immersion.
 Let $P$ be the intersection of $D_q$ with the closure
of $\beta\circ \alpha (Z)$ in $\bH_q$, i.e., $\dim P \leq m'+m''-1$. 
Since the natural $K$-action on $H_q$ extends to the action on $\bH_q$
so that its restriction to $D_q$ is transitive, $P$ does not meet $W_q$  for  a
general $(\beta, \alpha) \in K\times \cA$ by \cite[Theorem 1.15]{AFKKZ}. 
Hence, $\varphi_q|_{\beta \circ \alpha (Z)} : \beta \circ \alpha (Z) \to X$ is proper 
by Proposition \ref{cri.p1} and we are done.
\eproof

\bcor\label{main2.c1} Let $X$ be isomorphic (as an algebraic variety)
 either to a special linear group  $\SL_n(\bk)$  or to 
a symplectic group $\Sp_{2n} (\bk)$ and $Z$ be
an affine algebraic variety such that
$\ED (Z)\leq \dim X$. Then there exists a closed embedding of $Z$ into $X$.

\ecor

\bproof Suppose that $G'$ is the $\G_a^m$-subgroup of $\SL_n(\bk)$
(in particular, it is a unipotent abelian subgroup of a maximal dimension by \cite{Ma45})
as in the proof of Corollary \ref{main.c2} and $G''$ is the subgroup
that consists of the transposes of elements of $G'$. 
Note that $G'\cap G''=e$
 (where $e$ is the identity element of $G$) and $\dim G'=\dim G'' \geq \frac{ \dim X} {4}$. 
Hence, $\dim Z \leq \dim G'+\dim G''$
since $\ED (Z)\leq \dim X$ and, thus, $\dim Z \leq \frac{\dim X -1}{2}$.
 Similarly, for $X\simeq \Sp_{2n}(\bk)$ the maximal dimension of a unipotent abelian
subgroup $G'$ is greater than $\frac{ \dim X} {4}$ by \cite{Ma45} (see also \cite{Law}).
Furthermore, $G'$ can be chosen so that in a root space decomposition
its Lie algebra is generated by subspaces with positive roots \cite[page 7]{Law}. 
Replacing these positive roots by the corresponding negative roots we
get the Lie algebra of a maximal unipotent abelian subgroup $G''$ such
that $\dim G''=\dim G'$ and $G'\cap G''=e$. 
Hence,   $\dim Z \leq \dim G'+\dim G''$ as before
and Theorem \ref{main2.t1} implies the desired conclusion.
\eproof

In a more general setting we have the following.

\bcor\label{main2.c2} Let $Z$ be an affine algebraic variety, $X$ be an algebraic variety of the form 
$\A^{n_0}\times G_{1} \times  G_2 \times \ldots   \times G_l$
where each $G_i$ is
either  $\SL_{n_i}(\bk)$  or  $\Sp_{2n_i} (\bk)$.
Suppose that $\varphi : X \to Y$ is a finite morphism into a normal variety $Y$, $\ED (Z)\leq \dim Y$ and 
$S$ is a closed subvariety of $Y$ containing $ Y_{\sing}$ such that $\dim Z < \codim_Y S$.
Then $Z$ admits 
a closed embedding into $Y$ with the image contained in $Y\setminus S$. 
\ecor

\bproof By Theorem \ref{pre.t3} it suffices to consider the case of $Y=X$.
Since $X$ is isomorphic as an algebraic variety to a linear algebraic group 
$G=\G_a^{n_0}\times G_{1} \times  G_2 \times \ldots   \times G_l$
Theorem \ref{main2.t1} implies that
it is enough to construct $\G_a^m$-subgroups $G'$ and $G''$ of $G$ such that
$G'\cap G''=e$
 and $\dim Z\leq \dim G' + \dim G''$.
The proof of Corollary \ref{main2.c1} implies that one can find similar subgroups
$G_i'$ and $G_i''$ 
in each factor $G_i$ of $G$ such that $\dim G_i'+\dim G_i''\geq \frac{\dim G_i}{2}$. 
Thus, letting $G_i'=\G_a^{n_0}\oplus \bigoplus_{i=1}^l G_i'$ and $G_i''= \bigoplus_{i=1}^l G_i''$
we see that $\dim Z\leq \dim G' + \dim G''$ since $\dim Z \leq \frac{\dim G-1}{2}$.
This yields  the desired conclusion.
\eproof

\brem\label{main2.r1} If $G$ is a simple Lie group whose Dynkin diagram differs from $A_n$ or $C_n$, 
then there is no unipotent abelian
subgroup of $G$ whose dimension is at least $\frac{\dim G-1}{4}$ \cite{Ma45}. 
Hence, for such groups and  a smooth $Z$
our method   is less effective  than the one in \cite{FvS21}.
\erem



\begin{thebibliography}{KaMi} 


\bibitem[AFKKZ]{AFKKZ} I.~V.~Arzhantsev, H.~Flenner, S.~Kaliman,
F.~Kutzschebauch, M.~Zaidenberg,
{\em Flexible varieties and automophism groups}.  Duke Math.\ J.\ {\bf 162} (2013), no. 4, 767--823.

 






 

\bibitem[BMS]{BMS} S. Bloch, M. Pavaman Murthy, L. Szpiro, {\em Zero cycles and the number of generators of an ideal}, {\bf 38}, 1989, Colloque en 
l\'honneur de Pierre Samuel (Orsay, 1987), pp. 51-74.










%


 
 \bibitem[FvS21]{FvS21} P. Feller, I. van Santen, {\em  Existence of embedding of smooth varieties into linear algebraic groups}, 
 J. of Alg. Geom. (to appear), arXiv:2007.16164.

\bibitem[FKZ]{FKZ} H.~Flenner, S.~Kaliman, and M.~Zaidenberg, {\em A Gromov-Winkelmann type theorem for flexible varieties},  J. Eur. Math. Soc. (JEMS) 
{\bf18} (2016), no. 11, 2483-2510. 
%




\bibitem[Fr]{Fre} G. ~Freudenburg, {\em Algebraic Theory of Locally Nilpotent Derivations} , Encyclopaedia of Mathematical Sciences,  Springer,
 Berlin-Heidelberg-New York, 2006.




\bibitem[Gro58]{Gro58}  A.~ Grothendieck, {\em Torsion homologique et sections rationnelles}, Anneaux de Chow et Applications,
S\`eminaire Claude Chevalley, 1958, expos\^e n. 5, 1-29.





\bibitem[Hol]{Hol} A. Holme,  {\em Embedding-obstruction for singular algebraic varieties in $\PP^N$}, Acta Math. {\bf 135} (1975), no. 3-4, 155-185.

 



\bibitem[Ka91]{Ka91} S.~ Kaliman, {\em Extensions of isomorphisms between affine algebraic subvarieties of $k^n$ to automorphisms of $k^n$}, 
Proc. Amer. Math. Soc. {\bf 113} (1991), no. 2, 325-334.




\bibitem[Ka20]{Ka20} S.~Kaliman, {\em Extensions of isomorphisms of subvarieties in flexible varieties}, Transform. Groups {\bf 25} (2020), no. 2, 517-575.

\bibitem[Ka21]{Ka21} S.~Kaliman, {\em Lines in affine toric varieties}, Israel J. of Mathematics, TBD (2022) 1-29,
DOI 10.1007/s11856-022-2332-4.


 



\bibitem[KaUd]{KaUd} S. ~Kaliman, D. ~Udumyan, {\em On automorphisms of flexible varieties},
Adv. Math. 396 (2022), Paper No. 108112, 43 pp. 14R10 (14L30).






 

\bibitem[Law]{Law} R. Lawther, {\em Maximal abelian sets of roots},
 Mem. Amer. Math. Soc. {\bf 250} (2017), no. 1192, vii+219 pp.





\bibitem[Ma45]{Ma45}  A. Malcev,  {\em
Commutative subalgebras of semi-simple Lie algebras}, (Russian) Bull. Acad. Sci. URSS. S\'er. Math.
 [Izvestia Akad. Nauk SSSR] {\bf 9} (1945), 291-300. 








%
%
\bibitem[Ra]{Ra} C.~P.~Ramanujam, {\em A note on automorphism groups of algebraic varieties}, Math.\ Ann.\
{\bf 156} (1964), 25--33.

 \bibitem[ReYo]{ReYo} Z.~ Reichstein, B.~ Youssin, {\em Equivariant resolution of points of indeterminacy},
 Proc. Amer. Math. Soc. {\bf 130} (2002), no. 8, 2183 -2187.








\bibitem[Su]{Su} H.~ Sumihiro, {\em Equivarieant completion}, J. Math. Kyota Univ., {\bf 14:1} (1974), 1-14.



\bibitem[Sr]{Sr} V. Srinivas, {\em  On the embedding dimension of an affine variety},  Math. Ann., {\bf 289} (1991), no. 1, 25-132.

\bibitem[Swan]{Swan} R. G. Swan, {\em A cancellation theorem for projective modules in the metastable range}, Invent. Math. {\bf 27} (1974), 23-43.









\end{thebibliography}
\end{document}